\documentclass[11pt]{amsart}

\usepackage[all]{xy}
\usepackage{mathrsfs}
\usepackage{amsmath,amsthm, amssymb, amscd, amsgen,amsxtra,amsfonts, amsbsy, stmaryrd}
\usepackage[all]{xy}
\usepackage{longtable}
\usepackage{verbatim}
\usepackage{tikz}

\usepackage{color}
\definecolor{shadecolor}{gray}{0.875}
\definecolor{dblue}{rgb}{0,0,.6}
\usepackage[colorlinks=true, linkcolor=dblue, citecolor=dblue, filecolor = dblue, menucolor = dblue, urlcolor = dblue]{hyperref}

\usepackage{etex}
\usepackage{calligra,mathrsfs}
\usepackage{mathtools}
\usepackage{extarrows}
\usepackage{ifthen}

\usepackage{graphicx}
\usepackage{caption,rotating}
\usepackage{color}
\usepackage{subfigure}
\usepackage{wasysym}

\usepackage{times}

\usetikzlibrary{positioning}
\tikzset{main node/.style={circle,fill=blue!20,draw,minimum size=0.8cm,inner sep=0pt},}

\theoremstyle{definition}
\newtheorem{Definition}{Definition}[section]

\newtheorem{Remark}[Definition]{Remark}

\theoremstyle{plain}
\newtheorem{Theorem}[Definition]{Theorem}

\newtheorem{Corollary}[Definition]{Corollary}

\newtheoremstyle{voiditstyle}{3pt}{3pt}{\itshape}{\parindent}%
{\bfseries}{.}{ }{\thmnote{#3}}%
\theoremstyle{voiditstyle}

\newtheoremstyle{voidromstyle}{3pt}{3pt}{\rm}{\parindent}%
{\bfseries}{.}{ }{\thmnote{#3}}%
\theoremstyle{voidromstyle}

\newcommand{\EXT}{\mathop{\mathscr{E}{\kern -2pt {xt}}}\nolimits}

\newcommand{\HOM}{\mathop{\mathscr{H}{\kern -3pt {om}}}\nolimits}

\newcommand{\Frob}{\mathrm{F}}

\title{Frobenius actions on Del Pezzo surfaces of degree 2}
\author{Olof Bergvall}%
\address{Division of Mathematics and Physics, M\"alardalen University, 721 23 V\"aster\aa s, Sweden}
\email{olof.bergvall@mdu.se}

\date{\today}

\begin{document}
\begin{abstract}
 We determine the number of Del Pezzo surfaces of degree 2 over finite fields of odd characteristic with specified action of the Frobenius endomorphism,
 i.e. we solve the "quantitative inverse Galois problem". As applications we determine the number of Del Pezzo surfaces of degree 2 with a given number of points
 and recover results of Banwait-Fit\'e-Loughran and Loughran-Trepalin.
\end{abstract}
\maketitle

\section{Introduction}

Let $X$ be a smooth, projective, geometrically irreducible and geometrically rational surface
over a finite field $\mathbb{F}_q$ with $q$ elements. Let $\Frob \in \mathrm{Gal}(\overline{\mathbb{F}}_q/\mathbb{F}_q)$
denote the Frobenius element. We then have
$$
|X(\mathbb{F}_q)| = q^2+\mathrm{Tr}\left( F, \mathrm{Pic}(X_{\overline{\mathbb{F}}_q}) \right) \cdot q+1
$$
where $X(\mathbb{F}_q)$ denotes the set of $\mathbb{F}_q$-points of $X$
and $\mathrm{Tr}\left( F, \mathrm{Pic}(X_{\overline{\mathbb{F}}_q}) \right)$ denotes the trace of the Frobenius element
on the Picard group
of $X_{\overline{\mathbb{F}}_q}=X \otimes \overline{\mathbb{F}}_q$.

Over an algebraically closed field, any rational surface is isomorphic to the blow-up
of a minimal rational surface (i.e. the projective plane or a Hirzebruch surface).
We shall refer to the class of $X_{\overline{\mathbb{F}}_q}$ in this classification
as the type of $X$. Suppose that $X\cong \mathrm{Bl_{Z}}S$, where $S$ is a minimal 
rational surface. The realization of $X$ as a blow-up equips $X$ with a collection of
exceptional curves $E_1, \ldots, E_n$. This extra structure is called a geometric marking of $X$ (see \cite{dolgachevortland}, Section V.2).
For a given type $\mathrm{t}$ we denote the moduli space of rational surfaces
of type $\mathrm{t}$ by $\mathcal{M}_{\mathrm{t}}$ and we denote the moduli space
of geometrically marked rational surfaces of type $\mathrm{t}$ by $\mathcal{M}_{\mathrm{t}}^{\mathrm{gm}}$.
There is a natural forgetful morphism $\mathcal{M}_{\mathrm{t}}^{\mathrm{gm}} \to \mathcal{M}_{\mathrm{t}}$
and we denote its Galois group by $G_{\mathrm{t}}$. It is the group changing geometric markings (it is important to note that far from all permutations of a given geometric marking can be realized through a blow-up).

Let $a=\mathrm{Tr}\left( F, \mathrm{Pic}(X_{\overline{\mathbb{F}}_q}) \right)$. 
The value of $a$ completely determines the number of $\mathbb{F}_q$-points of $X$.
For a given type of rational surfaces, $a$ can only take a finite set of values.
It is thus very natural to ask:
\begin{enumerate}
\item For a given type of rational surfaces, which values of $a$ can occur?
\item For a given value of $a$, how many rational surfaces are there with $q^2+aq+1$ elements?
\end{enumerate}
Question (1) was first asked by Serre for the case of cubic surfaces, see \cite{serre}, Section 2.3.3.
Question (2) is a quantitative or refined version of Serre's question.

Suppose that $X$ is of type $\mathrm{t}$ and that the corresponding Galois group is $\mathrm{G}_{\mathrm{t}}$. The action of $\mathrm{Gal}(\overline{\mathbb{F}}_q/\mathbb{F}_q)$
on $X$ gives rise to a homomorphism $\mathrm{Gal}(\overline{\mathbb{F}}_q/\mathbb{F}_q) \to G_{\mathrm{t}}$ and thus to a Frobenius
element $\Frob \in G_{\mathrm{t}}$. Two natural refinements of the above two questions can now be formulated:
\begin{enumerate}
\setcounter{enumi}{2}
\item For a given type of rational surfaces, which conjugacy classes of $G$ can arise as the conjugacy class of $\Frob$?
\item For a given element $g$ of $G_{\mathrm{t}}$, how many rational surfaces are there such that $F$ is realized by $g$?
\end{enumerate}
Question (3) is sometimes called the "inverse Galois problem (of surfaces of given type)", see e.g. \cite{elsenhansjahnel}, Introduction, and \cite{loughrantrepalin}, Section 1.
Question (4) can thus be seen as a quantitative or refined version of the inverse Galois problem. Question (4) is clearly
the most ambitious one as its answer implies the answers to the other three questions.

Below, we recall what is known around the above questions for Del Pezzo surfaces of low degree (i.e. when there is moduli).

\subsection{Del Pezzo surfaces of degree 4}
For quartic Del Pezzo surfaces, Question (1) was solved by Banwait, Fit\'e and Loughran, see \cite{banwaitetal}, Theorem 1.3. Question (3) was solved by Trepalin, see \cite{trepalin} Theorem 1.4. Questions (2) and (4)
have been solved by Gounelas and the author, see \cite{bergvallgounelas}, Table 1. Thus, in the 
case of quartic Del Pezzo surfaces, the answers to Question (1)-(4) are known.

\subsection{Del Pezzo surfaces of degree 3}
Also in the case of cubic surfaces, the answers to Questions (1)-(4) are known.
Over $\mathbb{Q}$, Elsenhans and Jahnel solved Question (3), see \cite{elsenhansjahnel}. 
Over finite fields, Question (1) was solved by Banwait, Fit\'e and Loughran, see \cite{banwaitetal} Theorem 1.1. Question (3) was solved by Loughran and Trepalin, see \cite{loughrantrepalin} Theorem 1.1.
Questions (2) and (4) were solved by Das, see \cite{das_arithmetic} Table 1 (see also \cite{bergvallgounelas}).

\subsection{Del Pezzo surfaces of degree 2}
For Del Pezzo surfaces of degree 2, Question (1) has been solved by 
Fit\'e and Loughran, see \cite{banwaitetal} Theorem 1.4. Question (3)
has been solved by Loughran and Trepalin, see \cite{loughrantrepalin} Theorem 1.2.
Questions (2) and (4) are rather closely related to previous results of the author. The purpose of the present paper is
to extract these answers explicitly.

\subsection{Del Pezzo surfaces of degree 1}
For Del Pezzo surfaces of degree 2, Question (1) has been solved by 
Fit\'e and Loughran, see \cite{banwaitetal} Theorem 1.5. 
Questions (2)-(4) seem to be mostly open with the exception that
Mangano has counted the number of Del Pezzo surfaces of degree $1$
for which the Frobenius acts as the identity (excluding very small $q$, see \cite{mangano} Proposition 6.0.7).

\subsection{The present paper}
The purpose of the present paper is to record the answer to Question (4) in the case
of Del Pezzo surfaces over finite fields of odd characteristic, see Theorem~\ref{quant_class_cor}.
In other words, we determine the number of Del Pezzo surfaces $X$ of degree $2$ over $\mathbb{F}_q$ such that the Frobenius endomorphism acts on $X$ as a given element of the Weyl group $W(E_7)$ where $q$ is a power of an odd prime number.
As a consequence, 
we also obtain the answer to Question (2) and we recover the answers to Questions (1) and (3), see Corollaries~\ref{qual_class_cor}, \ref{quant_trace_cor} and \ref{qual_trace_cor}.

\begin{Remark}
We remark that the answers in characteristic 2 will not be exactly the same as those in odd characteristic
(e.g. a simple computation shows that there are no Del Pezzo surfaces of degree $2$ over $\mathbb{F}_2$
such that $\mathrm{Tr}\left( F, \mathrm{Pic}(X_{\overline{\mathbb{F}}_2})\right)=8$ but plugging in $q=2$ in our formula gives the result $135$).
This is rather expected given the results in \cite{glynn}, \cite{iss} and \cite{kklpw}.
\end{Remark}

\section{Del Pezzo surfaces of degree $2$}
Let $q$ be a power of an odd prime number $p$ and let $\mathbb{F_q}$ be a finite
field with $q$ elements. Let $\overline{\mathbb{F}}_q$ be an algebraic closure of
$\mathbb{F}_q$ and let $X$ be a smooth and complete surface of $\overline{\mathbb{F}}_q$.
The surface $X$ is called a Del Pezzo surface if the anticanonical class
$-K_X$ is ample. The number $d=K_X^2$ is called the degree of $X$; it is an integer
$1 \leq d \leq 9$. A Del Pezzo surface of degree $d$ is isomorphic to the blow-up of $\mathbb{P}^2$
in $n=9-d$ points (with the exception that $\mathbb{P}^1 \times \mathbb{P}^1$ cannot be obtained as a blow-up) in general position.

Thus, a Del Pezzo surface of degree $2$ can be realized as a blow-up of $\mathbb{P}^2$
in $7$ points $P_1, \ldots, P_7$ in general position, i.e. no three should lie on a line and
no six should lie on a conic. The choice of a representation of $X$ as a blow-up in
$7$ ordered points equips $X$ with $7$ exceptional classes $E_1, \ldots, E_7$.
Such a collection of classes (coming from a blow-up) is called a geometric marking, see \cite{dolgachevortland}.
If we denote the
class of a line in $\mathbb{P}^2$ by $L$ we have
$$
\mathrm{Pic}(X) = \mathbb{Z}L \oplus \mathbb{Z}E_1 \oplus \cdots \oplus \mathbb{Z}E_7
$$
with intersection pairing
$$
L^2 =1, \quad E_i^2=-1, \quad L \cdot E_i =0, \quad E_i \cdot E_j =0, \quad (i \neq j). 
$$
Any automorphism of the Picard group leaves the canonical class invariant. The orthogonal complement $K_X^{\perp}$
of the canonical class can be identified with the defining representation of the Weyl group $W(E_7)$.
More precisely, the classes
\begin{align*}
 & E_i-E_j, \quad 1 \leq i <j \leq 7 \\
 & L-E_i-E_j-E_k, \quad 1 \leq i <j<k \leq 7 \\
 & 2L - E_1- \cdots - E_7 + E_i, \quad 1 \leq i \leq 7
\end{align*}
constitute a choice of positive roots for the root system $E_7$
(the three types of classes above can naturally be interpreted in terms of the conditions that the
blown-up points should be distinct with no three on a line and no six on a conic).
We thus have a natural isomorphism
$$
\mathrm{Pic}(x) \cong V_{\mathrm{triv}} \oplus V_{\mathrm{std}}
$$
of $W(E_7)$-representations where $V_{\mathrm{triv}}$ denotes the trivial representation and
$V_{\mathrm{std}}$ denotes the standard representation of $W(E_7)$. 

A Del Pezzo surface $X$ of degree $2$ can also be realized as a double cover of
$\mathbb{P}^2$ branched over a smooth quartic curve $C$.
The image of a geometric marking $E_1, \ldots, E_7$ on $X$ is a collection
of $7$ ordered bitangents of $C$. Such $7$-tuples are known as Aronhold heptads
and they correspond to symplectic level $2$ structures on $C$ (see e.g. \cite{bergvall_arith}, \cite{bergvall_thesis} or \cite{dolgachevortland}). We denote the moduli space of
geometrically marked Del Pezzo surfaces by $\mathcal{M}_{\mathrm{gm}}$, we denote the moduli space of
plane quartics with symplectic level $2$ structure by $\mathcal{Q}[2]$ and
we denote the moduli space of ordered $7$-tuples of points in general position in
$\mathbb{P}^2$ by $\mathcal{P}^2_7$. There are then equivariant isomorphisms of coarse moduli spaces (see \cite{dolgachevortland}, Chapter IX, Theorem 1)
\begin{center}
\begin{tikzpicture}
    \node (M) at (0,0) {$\mathcal{M}^{\mathrm{gm}}$};
    \node (Q) at (-2,-2) {$\mathcal{Q}[2]$};
    \node (P) at (2,-2) {$\mathcal{P}^2_7$};
    \draw [stealth-stealth] (M) -- (Q);
    \draw [stealth-stealth] (Q) -- (P);
    \draw [stealth-stealth] (M) -- (P);
\end{tikzpicture} 
\end{center}
On the level of stacks, the bottom arrows are still isomorphisms while the others are not.
The reason is that every geometrically marked Del Pezzo surface of degree $2$ has precisely one non-trivial
automorphism; when given as a cover branched over a plane quartic, this automorphism can be seen as the involution switching the two sheets of the cover. Seven points in general position (up to projective equivalence) do not have any automorphisms and neither does a smooth quartic with symplectic level $2$ structure. 

The covering involution generates a subgroup $\mathbb{Z}/2\mathbb{Z}$ of $W(E_7)$ and we have an isomorphism $W(E_7) \cong \mathrm{Sp}(6,2) \times \mathbb{Z}/2\mathbb{Z}$
where $\mathrm{Sp}(6,2)$ denotes the symplectic group on a $6$-dimensional vector space over $\mathbb{F}_2$. The nontrivial element of $\mathbb{Z}/2\mathbb{Z}$
acts as $-1$ on $V_{\mathrm{std}}$. We therefore denote elements of the form $(g,1) \in W(E_7)$ simply by $g$ and we write $-g$ for elements of the form $(g,-1)$.
We use similar notation for conjugacy classes; we follow \cite{conwayetal}, p. 46-47, and denote a conjugacy class of $\mathrm{Sp}(6,2)$ by the order of its elements
followed by a letter used to distinguish conjugacy classes of the same order and we write a minus sign in front to denote the same class multiplied by $-1$. Thus,
$12A$ denotes a conjugacy class of $\mathrm{Sp}(6,2)$ consisting of elements of order $12$ and $-12A$ denotes the conjugacy class consisting of elements of type
$(g,-1)$ where $g \in 12A$.

\section{Lefschetz trace formula}

Let $Y$ be a scheme defined over $\mathbb{F}_q$ and let $G$ be a group of $\mathbb{F}_q$-automorphisms of $Y$.
Let $\ell$ be a prime number different from the characteristic and let $\mathrm{H}_{\text{\'et},c}^k\left( Y_{\overline{\mathbb{F}}_q},\mathbb{Q}_{\ell} \right)$
the $k$th compactly supported \'etale cohomology group of $Y_{\overline{\mathbb{F}}_q}$ with $\mathbb{Q}_{\ell}$ coefficients.
For an element $g \in G$, write $Y_{\overline{\mathbb{F}}_q}^g$ for the set of points of $Y_{\overline{\mathbb{F}}_q}$ fixed by $g$.
Lefschetz trace formula tells us that
$$
\left| Y^{\Frob g} \right| = \sum_{k \geq 0} (-1)^k \mathrm{Tr}\left(\Frob g, \mathrm{H}_{\text{\'et},c}^k\left( Y_{\overline{\mathbb{F}}_q},\mathbb{Q}_{\ell} \right) \right)
$$

We are trying to count the number of Del Pezzo surfaces $X$ of degree $2$ over $\mathbb{F}_q$ such that the Frobenius acts on $\mathrm{Pic}(X_{\overline{\mathbb{F}}_q})$ as $w$ for each element $w \in W(E_7)$. If $\Frob$ acts as $w$, then $\Frob w^{-1}$ acts as the identity and $X$ will then contribute to the count of $|\mathcal{M}_{\mathrm{gm}}^{\Frob w^{-1}}|$ and we will count $X$ once for every possible geometric marking on $X$, i.e. $|W(E_7)|$ times.

\section{Counts of surfaces}

In the previous section we saw that the problem of counting the number of Del Pezzo surfaces $X$ of
degree $2$ such that the Frobenius acts on $\mathrm{Pic}(X_{\overline{\mathbb{F}}_q})$ as $w \in W(E_7)$ is
intimately connected to the problem of finding the structure of the cohomology groups of the moduli space
$\mathcal{M}_{\mathrm{gm}}$ of geometrically marked Del Pezzo surfaces of degree $2$ as $W(E_7)$-representations.
The structure of the groups $\mathrm{H}_{\text{\'et},c}^k\left( \mathcal{M}_{\mathrm{gm}},\mathbb{Q}_{\ell} \right)$ as  $\mathrm{Sp}(6,2)$-representations
has been determined by the author in a series of papers, see \cite{BergvallGD}, \cite{BergvallEJM}, \cite{bergvalltor} and also \cite{bergvall_seven}.
However, we have noted that $W(E_7) \cong \mathrm{Sp}(6,2) \times \mathbb{Z}/2\mathbb{Z}$ and that the second factor acts trivially
on $\mathcal{M}_{\mathrm{gm}}$. Thus, the structure of $\mathrm{H}_{\text{\'et},c}^k\left( \mathcal{M}_{\mathrm{gm}},\mathbb{Q}_{\ell} \right)$ as a representation
of $W(E_7)$ is entirely determined by its structure as a representation of $\mathrm{Sp}(6,2)$ (in terms of characters, its values on a conjugacy class $c$ is the same as on the conjugacy class $-c$). We may thus use the results of \cite{BergvallGD}, Table 9, and Lefschetz trace formula to obtain the following.

\begin{Theorem}
\label{quant_class_cor}
 Let $\mathbb{F}_q$ be a finite field of odd characteristic and let $w$ be an element of
 $W(E_7)$. The number geometrically marked Del Pezzo surfaces of degree $2$ over $\mathbb{F}_q$ such that $w$ acts on
 $\mathrm{Pic}(X_{\overline{\mathbb{F}}_q})$ as $w$ is given 
 in Table~\ref{frob_action_tab}.
\end{Theorem}

As an immediate application we recover a result of Loughran and Trepalin 
(since we have assumed the characteristic to be odd we do not recover their results concerning characteristic $2$).

\begin{Corollary}[Loughran-Trepalin, \cite{loughrantrepalin} Theorem 1.2]
\label{qual_class_cor}
 Let $\mathbb{F}_q$ be a finite field of odd characteristic and let $c$ be a conjugacy class of
 $W(E_7)$. There is then a Del Pezzo surface $X$ of degree $2$ over $\mathbb{F}_q$ such that
 the Frobenius acts on $X$ as an element in $c$ except if
 \begin{enumerate}
  \item $c = \pm 1A$ and $q=3$, $5$ or $7$,
  \item $c = \pm 2A$ and $q=3$, or
  \item $c = \pm 2B$ and $q=3$ or $5$.
 \end{enumerate}
\end{Corollary}

Recall that the Picard group of a Del Pezzo surface $X$ of degree $2$ is isomorphic to
$$
\mathrm{Pic}(X_{\overline{\mathbb{F}}_q})  \cong V_{\mathrm{triv}} \oplus V_{\mathrm{std}}
$$
as a $W(E_7)$ representation. Thus, $a=\mathrm{Tr}\left( F, \mathrm{Pic}(X_{\overline{\mathbb{F}}_q}) \right)$
can be written as $a=1+b$ where $b$ is the trace of $F$ on $V_{\mathrm{std}}$.
Since $\mathrm{dim}(V_{\mathrm{std}})=7$, $b$ is an integer in the range $-7 \leq b \leq 7$
(not all these values occur though, see \cite{conwayetal} p. 47). Also recall that $W(E_7) \cong \mathrm{Sp}(6,2) \times
\mathbb{Z}/2\mathbb{Z}$. On $V_{\mathrm{triv}}$, the action is trivial
and on $V_{\mathrm{std}}$ the second factor acts as $-1$. 
Thus, if $g \in \mathrm{Sp}(6,2)$ has trace $a=1+b$ on $\mathrm{Pic}(X_{\overline{\mathbb{F}}_q})$,
then the trace of $-g$ is $1-b$. 

\begin{Corollary}
\label{quant_trace_cor}
 Let $\mathbb{F}_q$ be a finite field of odd characteristic. The number of Del Pezzo surfaces $X$ 
 of degree $2$ over $\mathbb{F}_q$ such that the Frobenius element acts with a given trace is as given in Table~\ref{trace_tab}.
\end{Corollary}

By considering the values of the character of $V_{\mathrm{triv}} \oplus V_{\mathrm{std}}$ we see that the possible
values of $a$ are
$$
-6, -4, -3, -2, -1, 0, 1, 2, 3, 4, 5, 6, 8
$$
We have therefore omitted the other values in Table~\ref{trace_tab}.
We observed above that if the trace of $g$ is $a=1+b$, then the trace of $-g$ is $1-b$; we thus expect the table to be symmetric around $a=1$ (i.e. there are as many Del Pezzo surfaces of degree $2$ such that
the Frobenius has trace $a$ as there are Del Pezzo surfaces of degree $2$ such that the Frobenius has trace $2-a$).

From Corollary~\ref{quant_trace_cor} we recover the following result
of Banwait, Fit\'e and Loughran (again, since we have assumed the characteristic to be odd we do not recover their results concerning characteristic $2$).

\begin{Corollary}[Banwait-Fit\'e-Loughran, \cite{banwaitetal} Theorem 1.4]
\label{qual_trace_cor}
 Let $\mathbb{F}_q$ be a finite field of odd characteristic and let
 $$
 a \in \{-6, -4, -3, -2, -1, 0, 1, 2, 3, 4, 5, 6, 8\}.
 $$
 There is then a Del Pezzo surface $X$ of degree $2$ over $\mathbb{F}_q$ such that
 the Frobenius acts on $X$ with trace $a$ except if
 \begin{enumerate}
 \item $a=-6$ or $a=8$ and $q=3$, $q=5$ or $q=7$, or
 \item $a=-4$ or $a=6$ and $q=3$.
 \end{enumerate}
 In particular, there are Del Pezzo surfaces over $\mathbb{F}_q$ such that
 the Frobenius acts with trace $a$ for all $-3 \leq a \leq 5$.
\end{Corollary}

\begin{center}
\begin{table}
\begin{tabular}{l|l} 
 Conjucagy class of $w$ & Number of surfaces \\
 \hline
 $\pm 1A$ & $(q^3 - 20q^2 + 119q - 175)(q - 3)(q - 5)(q - 7)$ \\
 $\pm 2A$ & $(q^4 - 11q^3 + 43q^2 - 60q + 5)(q - 1)(q - 3)$ \\
 $\pm 2B$ & $(q^4 - 3q^3 - 13q^2 + 16q + 41)(q - 3)(q - 5)$ \\
 $\pm 2C$ & $(q^3 - 4q^2 - q + 5)(q + 1)(q - 1)(q - 3)$ \\
 $\pm 2D$ & $(q^5 - 2q^4 - 8q^3 + 11q^2 + 17q - 7)(q - 1)$ \\
 $\pm 3A$ & $(q^2 - 3q + 5)(q + 1)(q - 1)(q - 2)q$ \\
 $\pm 3B$ & $(q^2 + 2q + 3)(q + 2)(q + 1)(q - 2)^2$ \\
 $\pm 3C$ & $q^6 - 2q^5 - 2q^4 - 8q^3 + 16q^2 + 10q + 21$ \\
 $\pm 4A$ & $(q^5 - 2q^3 - 3q^2 - 3q + 3)(q + 1)$ \\
 $\pm 4B$ & $(q^3 - 4q^2 + 3q + 3)(q + 1)(q - 1)^2$ \\
 $\pm 4C$ & $(q^4 - q^3 - q^2 - 3)(q + 1)(q - 1)$ \\
 $\pm 4D$ & $q^6 - 3q^5 - 2q^4 + 7q^3 - 6q^2 + 16q + 11$ \\
 $\pm 4E$ & $(q^4 - q^3 - q^2 - 3)(q + 1)(q - 1)$ \\
 $\pm 5A$ & $(q^2 + 1)(q + 1)(q - 1)q^2$ \\
 $\pm 6A$ & $(q^2 - q - 1)(q + 1)(q - 1)(q - 2)q$ \\
 $\pm 6B$ & $(q^2 - q - 1)(q + 2)(q + 1)(q - 1)q$ \\
 $\pm 6C$ & $(q^4 - q^2 - 2q - 4)(q + 1)q$ \\
 $\pm 6D$ & $(q^3 - q^2 - q + 2)(q + 1)(q - 1)q$ \\
 $\pm 6E$ & $(q^5 + q^4 + q^3 - 3q^2 - q - 3)(q - 1)$ \\
 $\pm 6F$ & $(q^5 - q^4 + q^3 - 5q^2 + q - 3)(q - 1)$ \\
 $\pm 6G$ & $(q^2 + q + 1)^2(q - 1)^2$ \\
 $\pm 7A$ & $(q^2 - q + 1)(q + 1)q^3$ \\
 $\pm 8A$ & $(q^3 - q + 1)(q^2 + 1)(q + 1)$ \\
 $\pm 8B$ & $(q^3 - q - 1)(q^2 + 1)(q - 1)$ \\
 $\pm 9A$ & $(q^2 + q + 1)(q^2 - q + 1)(q + 1)q$ \\
 $\pm 10A$ & $(q^2 + 1)(q + 1)(q - 1)q^2$ \\
 $\pm 12A$ & $(q^2 + q + 1)(q + 1)(q - 1)q^2$ \\
 $\pm 12B$ & $(q^2 - q - 1)(q + 1)(q - 1)q^2$ \\
 $\pm 12C$ & $(q^2 + 1)(q + 1)q^3$ \\
 $\pm 15A$ & $(q^2 + 1)(q + 1)(q - 1)q^2$
\end{tabular}
\caption{The number of geometrically marked Del Pezzo surfaces $X$ of degree $2$ over $\mathbb{F}_q$ such that the Frobenius element
acts on $\mathrm{Pic}(X_{\overline{\mathbb{F}}_q})$ as $w \in W(E_7)$, see \cite{conwayetal} p. 46-47 for notation.}
\end{table}
\label{frob_action_tab}
\end{center}

\begin{center}
\begin{table}
\begin{tabular}{r|l}
Trace & Number of surfaces \\
\hline
 $-6$ & $(q^3 - 20q^2 + 119q - 175)(q - 3)(q - 5)(q - 7)$ \\
 $-4$ & $63(q^4 - 11q^3 + 43q^2 - 60q + 5)(q - 1)(q - 3)$ \\
 $-3$ & $672(q^2 - 3q + 5)(q + 1)(q - 1)(q - 2)q$ \\
 $-2$ & $945(13q^5 - 24q^4 + 10q^3 - 7q^2 - 59q + 51)(q + 1)$ \\
 $-1$ & $896(169q^5 - 99q^4 + 19q^3 - 159q^2 - 80q + 60)(q + 1)$ \\
 $0$ & $722883q^6 - 78225q^5 + 132510q^4 - 333375q^3 + 301602q^2 + 116760q + 1330245$ \\
 $1$ & $1728(653q^4 - 618q^3 + 548q^2 + 7q - 70)(q + 1)q$ \\
 $2$ & $722883q^6 - 78225q^5 + 132510q^4 - 333375q^3 + 301602q^2 + 116760q + 1330245$ \\
 $3$ & $896(169q^5 - 99q^4 + 19q^3 - 159q^2 - 80q + 60)(q + 1)$ \\
 $4$ & $945(13q^5 - 24q^4 + 10q^3 - 7q^2 - 59q + 51)(q + 1)$ \\
 $5$ & $672(q^2 - 3q + 5)(q + 1)(q - 1)(q - 2)q$ \\
 $6$ & $63(q^4 - 11q^3 + 43q^2 - 60q + 5)(q - 1)(q - 3)$ \\
 $8$ & $(q^3 - 20q^2 + 119q - 175)(q - 3)(q - 5)(q - 7)$
\end{tabular}
\caption{The number of Del Pezzo surfaces $X$ of degree $2$ over $\mathbb{F}_q$ such that the Frobenius element
acts on $\mathrm{Pic}(X_{\overline{\mathbb{F}}_q})$ with trace $a$ for all possible values of $a$.}
\label{trace_tab}
\end{table}
\end{center}

\clearpage

\bibliographystyle{acm}

\renewcommand{\bibname}{References}

\bibliography{references}

\end{document}